\documentclass[11pt]{article}

\usepackage{amsmath,amssymb}
\usepackage[mathscr]{eucal}
\usepackage{theorem}

\makeatletter
\@addtoreset{equation}{section}
\makeatother

\newtheorem{theorem}{Theorem}[section]
\newtheorem{proposition}[theorem]{Proposition}
\newtheorem{lemma}[theorem]{Lemma}
\newtheorem{corollary}[theorem]{Corollary}
{\theorembodyfont{\rmfamily}
\newtheorem{example}[theorem]{Example}
\newtheorem{remark}[theorem]{Remark}}

\title{\Huge\bf A q-Analog of the Hua Equations}
\author{\Large O. Bershtein\thanks{Partially supported by a grant of the
President of Ukraine} \and \Large S. Sinel'shchikov}
\date{Mathematics Division, B. Verkin Institute for Low Temperature Physics
and Engineering,
\\ National Academy of Sciences of Ukraine
\\ 47 Lenin Ave., Kharkov 61103, Ukraine}

\begin{document}
\large

\maketitle

\begin{center}
E-mail: bershtein@ilt.kharkov.ua
\\ \hspace{5em} sinelshchikov@ilt.kharkov.ua
\end{center}

\begin{abstract}
A necessary condition is established for a function to be in the image of a
quantum Poisson integral operator associated to the Shilov boundary of the
quantum matrix ball. A quantum analogue of the Hua equations is introduced.

\vspace{3ex} {\it Key words}: quantum matrix ball, Shilov boundary, Poisson
integral operator, invariant kernel, Hua equations.

(secondary).
\end{abstract}

\hfill {\sl Dedicated to the memory of L. L. Vaksman}

\bigskip

\section{Introduction}

In late 90ies three groups of specialists advanced in putting the basics of
the quantum theory of bounded symmetric domains.

T. Tanisaki and his team introduced $q$-analogs for prehomogeneous vector
spaces of commutative parabolic type and found an explicit form of the
associated Sato-Bernstein polynomials \cite{KMT, Kamita0, Morita, Kamita3}.

On the other hand, H. Jakobsen suggested a less intricate method of
producing the above quantum vector spaces. Actually he was on the way to
quantum Hermitian symmetric spaces of non-compact type \cite{Jak-Hermit,
JakZhang}. A similar approach is used in the work by W. Baldoni and P.
Frajria \cite{Baldoni_Frajria} on $q$-analogs of algebras of invariant
differential operators and the Harish-Chandra homomorphism for these
quantum symmetric spaces. During the same period Jakobsen obtained a
description of all the unitarizable highest weight modules over the
Drinfeld-Jimbo algebras \cite{Jak3}.

The authors listed above made no use of the full symmetry of the quantum
prehomogeneous vector spaces in question \cite{SV0, ShSinVak2001}, which
became an obstacle in producing the quantum theory of bounded symmetric
domains.

The paper \cite{SV} laid the foundations of this theory. The subsequent
results have been obtained in the works of L. Vaksman, D. Proskurin, S.
Sinel'shchikov, A. Stolin, D. Shklyarov, L. Turowska, H. Zhang \cite{Turow,
ProskTurow, ShklZhang, VakShkl, SSV_LMP, Vak01, SSV01}. The compatibility
of the approaches described above \cite{Tanisaki, Jak3, SV} has been proved
by Shklyarov \cite{Shkl}. The study of quantum analogs for the
Harish-Chandra modules related to quantum bounded symmetric domains and
their geometric realizations has been started in \cite{SSSV1, SStV, SSSV2}.
The present work proceeds with this research.

Recall that a bounded domain $\mathbb{D}$ in a finite dimensional vector
space is said to be symmetric if every point $p\in\mathbb{D}$ is an
isolated fixed point of a biholomorphic involutive automorphism
$\varphi_p:\mathbb{D}\to\mathbb{D}$,
$\varphi_p\circ\varphi_p=\operatorname{id}$.

Equip the vector space of linear maps in $\mathbb{C}^n$ and the canonically
isomorphic vector space $\operatorname{Mat}_n$ of complex $n\times n$
matrices with the operator norms. It is known \cite{Helg} that the unit
ball $\mathbb{D}=
\{\mathbf{z}\in\operatorname{Mat}_n|\:\mathbf{z}\mathbf{z}^*<1\}$ is a
bounded symmetric domain.

 Denote by $S(\mathbb{D})$ the Shilov boundary of $\mathbb{D}$,
$S(\mathbb{D})=\{\mathbf{z}\in\operatorname{Mat}_n|\: I-\mathbf{zz}^*=0\}\cong U_n$. It
is well known that both $\mathbb{D}$ and $S(\mathbb{D})$ are homogeneous spaces of the
group $SU_{n,n}$. Consider a function on $\mathbb{D}\times S(\mathbb{D})$ given by
\begin{equation}\label{classical_Poisson}
P(\mathbf{z},\mathbf{\zeta})=
\frac{\det(1-\mathbf{zz}^*)^n}{|\det(1-{\mathbf{\zeta}^*\mathbf{z}})|^{2n}},
\qquad\mathbf{z}\in\mathbb{D},\quad\mathbf{\zeta}\in S(\mathbb{D}).
\end{equation}
It is called the Poisson kernel \cite{Hua} associated to the Shilov
boundary.

The general concepts on boundaries of Hermitian symmetric spaces of
non-compact type and the associated Poisson kernels are exposed in
\cite{Koranyi_Inv, Koranyi_harmonic}.

The Poisson kernel, together with the $S(U_n\times U_n)$-invariant integral
$\int\limits_{U_n}\cdot d\nu(\zeta)$ on $U_n$, allow one to define the
Poisson integral operator
\begin{equation*}
\mathcal{P}:f(\zeta)\mapsto\int\limits_{U_n}P(\mathbf{z},\zeta)f(\zeta)
d\nu(\zeta),\qquad\mathbf{z}\in\mathbb{D},\quad\zeta\in S(\mathbb{D}).
\end{equation*}
It intertwines the actions of $SU_{n,n}$ in the spaces of continuous
functions on the domain $\mathbb{D}$ and on the Shilov boundary
$S(\mathbb{D})$. However, not every continuous function on $\mathbb{D}$ can
be produced via applying the Poisson integral operator to a continuous
function on the Shilov boundary $S(\mathbb{D})$.

Hua \cite{Hua} obtained the initial results on differential equations whose
solutions include the functions of the form $\mathcal{P}(f)$. A later
result of Johnson and Kor\'anyi \cite{JohnsonKoranyi} provided a system of
differential equations which give a complete characterization of such
functions. Their version of the Hua equations is as follows:
\begin{eqnarray*}\label{classical_Hua_1}
\sum_{i,j,k=1}^n\left(\delta_{ij}-\sum_{c=1}^nz_i^c\overline{z}_j^c\right)
\left(\delta_{k\alpha}-\sum_{c=1}^n\overline{z}_c^kz_c^\alpha\right)
\frac{\partial^2}{\partial\overline{z}_j^k\partial z_i^a}
u(\mathbf{z}) &=& 0,
\\ \label{classical_Hua_2}
\sum_{i,j,k=1}^n\left(\delta_{ai}-\sum_{c=1}^nz_a^c\overline{z}_i^c\right)
\left(\delta_{jk}-\sum_{c=1}^n\overline{z}_c^jz_c^k\right)
\frac{\partial^2}{\partial\overline{z}_i^j\partial z_\alpha^k}
u(\mathbf z) &=& 0
\end{eqnarray*}
for $a,\alpha=1,\ldots,n$.

This can be also represented in the form:
\begin{equation}\label{classical_Hua_new1}
\left.\sum_{c=1}^n\frac{\partial^2u(g\cdot\mathbf{z})}
{\partial z_c^\beta\partial\overline{z}_c^\alpha}
\right|_{\mathbf{z}=0}=0,\qquad g\in SU_{n,n},\quad
\alpha,\beta\in\{1,2,\ldots,n\}.
\end{equation}

It is known that the Poisson kernel \eqref{classical_Poisson} as a function
of $\mathbf{z}$ is a solution of the above equation system. While proving
this, one may stick to the special case $g=1$ because
$P(g\mathbf{z},\zeta)=\mathrm{const}(g,\zeta)P(\mathbf{z},g^{-1}\zeta)$,
cf. \cite[p. 597]{JohnsonKoranyi}. What remains is to note that the
relations
\begin{equation}\label{Hua_zero_classical1}
\left.\sum_{c=1}^n\frac{\partial^2P}
{\partial z_c^\beta\partial\overline{z}_c^\alpha}
\right|_{\mathbf{z}=0}=0,\qquad\zeta\in U_n,\quad
\alpha,\beta\in\{1,2,\ldots,n\}
\end{equation}
follow from \eqref{classical_Poisson}.

\medskip

Thus in the classical case the Poisson integral operator applied to a
function on the Shilov boundary is a solution of
\eqref{classical_Hua_new1}. We are going to obtain a quantum analog of this
well known result.

This work suggests a quantum analog for the Hua equations. We thus get a
quantization of the necessary condition for a function being a Poisson
integral of a function on the Shilov boundary.

The authors acknowledge that the ideas of this work have been inspired by
plentiful long term communications with L. Vaksman. Also, we would like to
express our gratitude to D. Shklyarov for pointing out several inconsistencies in a
previous version of this work, and encouraging us to improve the exposition.

\bigskip

\section{A background on function theory in quantum matrix
ball}\label{prelim}

In what follows we assume $\mathbb{C}$ to be the ground field and all
the algebras are assumed associative and unital.

Recall the construction of the quantum universal enveloping algebra for the
Lie algebra $\mathfrak{sl}_N$. The quantum universal enveloping algebras
were introduced by V. Drinfeld and M. Jimbo in an essentially more general
way than it is described below. We follow the notation of \cite{BrGood,
DaDeCon, Jant, Rosso}.

Let $q\in(0,1)$. The Hopf algebra $U_q\mathfrak{sl}_{N}$ is given by its
generators $K_i$, $K_i^{-1}$, $E_i$, $F_i$, $i=1,2,\ldots,N-1$, and the
relations:
$$K_iK_j=K_jK_i,\qquad K_iK_i^{-1}=K_i^{-1}K_i=1,$$
$$K_iE_j=q^{a_{ij}}E_jK_i,\qquad K_iF_j=q^{-a_{ij}}F_jK_i,$$
$$E_iF_j-F_jE_i=\delta_{ij}\frac{K_i-K_i^{-1}}{q-q^{-1}},$$
$$E_i^2E_j-(q+q^{-1})E_iE_jE_i+E_jE_i^2=0,\qquad |i-j|=1,$$
$$F_i^2F_j-(q+q^{-1})F_iF_jF_i+F_jF_i^2=0,\qquad |i-j|=1,$$
$$E_iE_j-E_jE_i=F_iF_j-F_jF_i=0,\qquad |i-j|\ne 1,$$
with $a_{ii}=2$, $a_{ij}=-1$ for $|i-j|=1$, $a_{ij}=0$ otherwise, and the
comultiplication $\Delta$, the antipode $S$, and the counit $\varepsilon$
being defined on the generators by
$$
\Delta(E_i)=E_i\otimes 1+K_i\otimes E_i,\quad
\Delta(F_i)=F_i\otimes K_i^{-1}+1\otimes F_i,\quad
\Delta(K_i)=K_i\otimes K_i,
$$
\begin{equation*}
S(E_i)=-K_i^{-1}E_i,\qquad S(F_i)=-F_iK_i,\qquad S(K_i)=K_i^{-1},
\end{equation*}
$$\varepsilon(E_i)=\varepsilon(F_i)=0,\qquad\varepsilon(K_i)=1,$$
see also \cite[Chapter 4]{Jant}.

Consider the Hopf algebra $U_q\mathfrak{sl}_{2n}$. Equip $U_q\mathfrak{sl}_{2n}$ with a
structure of Hopf $*$-algebra determined by the involution
$$K_j^*=K_j,\qquad E_j^*=
\begin{cases}
K_jF_j, & j\ne n,
\\ -K_jF_j, & j=n,
\end{cases}\qquad
F_j^*=
\begin{cases}
E_jK_j^{-1}, & j\ne n,
\\ -E_jK_j^{-1}, & j=n.
\end{cases}
$$
This Hopf $*$-algebra $(U_q\mathfrak{sl}_{2n},*)$ is denoted by
$U_q\mathfrak{su}_{n,n}$.

Denote by
$U_q\mathfrak{k}=U_q\mathfrak{s}(\mathfrak{u}_n\times\mathfrak{u}_n)$
the Hopf $*$-subalgebra generated by
$$E_i,F_i,\quad i\ne n;\qquad K_j^{\pm 1},\quad j =1,2,\ldots,2n-1.$$

\bigskip

Now we introduce the notation to be used in the sequel, and recall some known results on
the quantum matrix ball $\mathbb{D}$.

Consider a $*$-algebra $\operatorname{Pol}(\operatorname{Mat}_n)_q$ with generators
$\{z_a^\alpha\}_{a,\alpha=1,2,\ldots,n}$ and defining relations
\begin{equation}\label{z}
z_a^\alpha z_b^\beta=
\begin{cases}
qz_b^\beta z_a^\alpha, & a=b\;\&\;\alpha<\beta\;\text{\ or\ }\;
a<b\;\&\;\alpha=\beta
\\ z_b^\beta z_a^\alpha, & a<b \;\&\;\alpha>\beta
\\ z_b^\beta z_a^\alpha+(q-q^{-1})z_a^\beta z_b^\alpha, & a<b\;\&\;
\alpha<\beta,
\end{cases}
\end{equation}
$$
(z_b^\beta)^*z_a^\alpha=q^2 \sum_{a',b'=1}^n\sum_{\alpha',\beta'=1}^m
R(b,a,b',a')R(\beta,\alpha,\beta',\alpha')
z_{a'}^{\alpha'}(z_{b'}^{\beta'})^*+(1-q^2)\delta_{ab} \delta^{\alpha\beta}
$$
with $\delta_{ab}$, $\delta^{\alpha\beta}$ being the Kronecker symbols and
\begin{equation*}
R(b,a,b',a')=
\begin{cases}
q^{-1}, & a\ne b\;\&\;b=b'\;\&\;a=a'
\\ 1, & a=b=a'=b'
\\ -(q^{-2}-1), & a=b\;\&\;a^\prime=b'\;\&\;a'>a
\\ 0, &\text{otherwise}.
\end{cases}
\end{equation*}

Denote by $\mathbb{C}[\operatorname{Mat}_n]_q$ the subalgebra generated by
$z_a^\alpha$, $a,\alpha=1,2,\ldots,n$. It is a very well known quantum
analog of the algebra of holomorphic polynomials on $\operatorname{Mat}_n$.

Consider an arbitrary Hopf algebra $A$ and an $A$-module algebra $F$.
Suppose that $A$ is a Hopf $*$-algebra. The $*$-algebra $F$ is said to be
an $A$-module algebra if the involutions are compatible as follows:
\begin{equation*}
(af)^*=(S(a))^*f^*,\qquad a\in A,\quad f\in F,
\end{equation*}
with $S$ being the antipode of $A$.

It was demonstrated in \cite{SSV4}, see Propositions 8.12 and 10.1, that one has

\begin{proposition}\label{z_gener}
$\mathbb{C}[\operatorname{Mat}_n]_q$ carries a structure of
$U_q\mathfrak{sl}_{2n}$-module algebra given by
$$
K_n^{\pm 1}z_a^\alpha=
\begin{cases}
q^{\pm 2}z_a^\alpha, & a=\alpha=n,
\\ q^{\mp 1}z_a^\alpha, & a=n\;\&\;\alpha\ne n\text{\ or\ }
a\ne n\;\&\;\alpha=n,
\\ z_a^\alpha, & \text{otherwise},
\end{cases}
$$
$$
F_nz_a^\alpha=q^{1/2}\cdot
\begin{cases}
1, & a=\alpha=n,
\\ 0, & \text{otherwise},
\end{cases}
$$
$$
E_nz_a^\alpha=-q^{1/2}\cdot
\begin{cases}
q^{-1}z_a^nz_n^\alpha, & a\ne n\;\&\;\alpha\ne n,
\\ (z_n^n)^2, & a=\alpha=n,
\\ z_n^nz_a^\alpha, & \text{otherwise},
\end{cases}
$$
and for $k\ne n$
$$
K_k^{\pm 1}z_a^\alpha=
\begin{cases}
q^{\pm 1}z_a^\alpha, & k<n\;\&\;a=k\text{\ or\ }k>n\;\&\;\alpha=2n-k,
\\ q^{\mp 1}z_a^\alpha, & k<n\;\&\;a=k+1\text{\ or\ }k>n\;\&\;\alpha=2n-k+1,
\\ z_a^\alpha, & \text{otherwise with\ }k\ne n;
\end{cases}
$$
$$
F_kz_a^\alpha=q^{1/2}\cdot
\begin{cases}
z_{a+1}^\alpha, & k<n\;\&\;a=k,
\\ z_a^{\alpha+1}, & k>n\;\&\;\alpha=2n-k,
\\ 0, & \text{otherwise with\ }k\ne n;
\end{cases}
$$
$$
E_kz_a^\alpha=q^{-1/2}\cdot
\begin{cases}
z_{a-1}^\alpha, & k<n\;\&\;a=k+1,
\\ z_a^{\alpha-1}, & k>n\;\&\;\alpha=2n-k+1,
\\ 0, & \text{otherwise with\ }k\ne n
\end{cases}
$$
Also $\mathrm{Pol}(\mathrm{Mat}_n)_q$ is equipped this way with a structure
of $U_q\mathfrak{su}_{n,n}$-module algebra.
\end{proposition}

It is well known that in the classical case $q=1$ the Shilov boundary of
the matrix ball $\mathbb{D}$ is just the set $S(\mathbb{D})$ of all unitary
matrices. Our intention is to produce a q-analogue of the Shilov boundary
for the quantum matrix ball. Introduce the notation for the quantum minors
of the matrix $\mathbf{z}=(z_a^\alpha)$:
$$
(z^{\wedge k})_{\{a_1,a_2,\ldots,a_k
\}}^{\{\alpha_1,\alpha_2,\ldots,\alpha_k
\}}\stackrel{\mathrm{def}}{=}\sum_{s \in
S_k}(-q)^{l(s)}z_{a_1}^{\alpha_{s(1)}}z_{a_2}^{\alpha_{s(2)}}\cdots
z_{a_k}^{\alpha_{s(k)}},
$$
with $\alpha_1<\alpha_2<\ldots<\alpha_k$, $a_1<a_2<\ldots<a_k$, and $l(s)$
being the number of inversions in $s\in S_k$.

It is well known that the quantum determinant
$$
\det \nolimits_q \mathbf{z}=(z^{\wedge n})_{\{1,2,\ldots,n
\}}^{\{1,2,\ldots,n \}}
$$
is in the center of $\mathbb{C}[\mathrm{Mat}_n]_q$. The localization
of $\mathbb{C}[\mathrm{Mat}_n]_q$ with respect to the multiplicative
system $(\det_q \mathbf{z})^\mathbb{N}$ is called the algebra of
regular functions on the quantum $GL_n$ and is denoted by
$\mathbb{C}[GL_n]_q$.

\begin{lemma}[Lemma 2.1 of \cite{Vak01}]\label{inv}
There exists a unique involution $*$ in $\mathbb{C}[GL_n]_q$ such
that
$$
(z_a^\alpha)^*=(-q)^{a+\alpha-2n}(\det \nolimits_q
\mathbf{z})^{-1}\det \nolimits_q \mathbf{z}_a^\alpha,
$$
with $\mathbf{z}_a^\alpha$ being the matrix derived from $\mathbf{z}$ via
deleting the row $\alpha$ and the column $a$.
\end{lemma}

The $*$-algebra $\mathbb C[S(\mathbb{D})]_q=(\mathbb{C}[GL_n]_q,*)$ is a
q-analogue of the algebra of regular functions on the Shilov boundary of
the matrix ball $\mathbb{D}$. It can be verified easily that $\mathbb
C[S(\mathbb{D})]_q$ is a $U_q\mathfrak{su}_{n,n}$-module algebra (see
\cite[Theorem 2.2 and Proposition 2.7]{Vak01} for the proof).


There exists another definition of the algebra $\mathbb
C[S(\mathbb{D})]_q$.

Consider the two-sided ideal $J$ of the $*$-algebra
$\operatorname{Pol}(\operatorname{Mat}_n)_q$ generated by the
relations
\begin{equation}\label{Shilov_eq}
\sum\limits_{j=1}^n q^{2n-\alpha-\beta}z^\alpha_j(z^\beta_j)^*
-\delta^{\alpha\beta}=0,\qquad\alpha,\beta=1,2,\ldots,n.
\end{equation}
One can prove that this ideal is $U_q\mathfrak{su}_{n,n}$-invariant, which allows to
introduce the $U_q\mathfrak{su}_{n,n}$-module algebra
$\operatorname{Pol}(\operatorname{Mat}_n)_q/J$. It was proved in \cite[p. 381,
Proposition 6.1]{Vak01} that $\mathbb
C[S(\mathbb{D})]_q=\operatorname{Pol}(\operatorname{Mat}_n)_q/J$. The last equality also
works as a definition for the algebra of regular functions on the Shilov boundary of the
quantum matrix ball.

\bigskip

A module $V$ over $U_q\mathfrak{sl}_{2n}$ is said to be a weight module if
\begin{equation*}\label{weight}
V=\bigoplus\limits_{\mathbf{\lambda}\in P}V_{\mathbf{\lambda}},\qquad
V_\lambda=\left\{v\in V\;\left|\;K_iv=q^{\lambda_i}v,\quad
i=1,2,\ldots,2n-1\right.\right\},
\end{equation*}
where $\mathbf{\lambda}=(\lambda_1,\lambda_2,\ldots,\lambda_{2n-1})$ and
$P\cong\mathbb{Z}^{2n-1}$ is the weight lattice of the Lie algebra
$\mathfrak{sl}_{2n}$. A non-zero summand $V_{\mathbf{\lambda}}$ in this
decomposition is called the weight subspace for the weight $\lambda$.

We associate to every weight $U_q\mathfrak{sl}_{2n}$-module $V$ linear maps
$H_i$, $i=1,2,\ldots,2n-1$, in $V$ such that
\begin{equation*}
H_iv=\lambda_iv,\qquad \text{iff}\qquad v\in V_\lambda.
\end{equation*}

Fix the element
\begin{equation*}
H_0=\sum\limits_{j=1}^{n-1}j(H_j+H_{2n-j})+nH_n.
\end{equation*}
Any weight $U_q\mathfrak{sl}_{2n}$-module $V$ can be equipped with a $\mathbb{Z}$-grading
$V=\mathop{\oplus}\limits_r V_r$ by setting $v\in V_r$ if $H_0v=2rv$.

In what follows some more sophisticated spaces will be used. It is known
that
$$
\operatorname{Pol}(\operatorname{Mat}_n)_q=\bigoplus\limits_{k,j=0}^\infty
\mathbb{C}[\operatorname{Mat}_n]_{q,k}\cdot
\mathbb{C}[\overline{\operatorname{Mat}}_n]_{q,-j}.
$$
Here $\mathbb{C}[\overline{\operatorname{Mat}}_n]_q$ is the subalgebra of
$\operatorname{Pol}(\mathrm{Mat}_n)_q$ generated by $(z_a^\alpha)^*$,
$a,\alpha=1,2,\ldots,n$, and
$\mathbb{C}[\overline{\operatorname{Mat}}_n]_{q,-j}$,
$\mathbb{C}[\operatorname{Mat}_n]_{q,k}$ are the homogeneous components
related to the grading
$$
\deg(z_a^\alpha)=1,\quad\deg(z_a^\alpha)^*=-1,\qquad a,\alpha=1,2,\ldots,n.
$$
To rephrase this, every $f\in\operatorname{Pol}(\operatorname{Mat}_n)_q$ is
uniquely decomposable as a finite sum
\begin{equation}\label{series_f}
f=\sum\limits_{k,j\ge 0}f_{k,j},\qquad f_{k,j}\in
\mathbb{C}[\operatorname{Mat}_n]_{q,k}\cdot
\mathbb{C}[\overline{\operatorname{Mat}}_n]_{q,-j}.
\end{equation}
Note that $\dim\mathbb{C}[\operatorname{Mat}_n]_{q,k}\cdot
\mathbb{C}[\overline{\operatorname{Mat}}_n]_{q,-j}<\infty$.

Consider the vector space $\mathscr{D}(\mathbb{D})_q'$ of formal series of
the form \eqref{series_f} with the termwise topology. The
$U_q\mathfrak{sl}_{2n}$-action and the involution $*$ admit an extension by
continuity from the dense linear subspace
$\operatorname{Pol}(\operatorname{Mat}_n)_q$ to
$\mathscr{D}(\mathbb{D})_q'$
\begin{equation*}
*:\sum\limits_{k,j=0}^\infty f_{k,j}\mapsto\sum\limits_{k,j=0}^\infty
f_{k,j}^*.
\end{equation*}
Moreover, $\mathscr{D}(\mathbb{D})_q'$ is a $U_q\mathfrak{sl}_{2n}$-module bimodule over
$\operatorname{Pol}(\operatorname{Mat}_n)_q$. We call the elements of
$\mathscr{D}(\mathbb{D})_q'$ distributions on a quantum bounded symmetric domain.

\bigskip

\section{Statement of the main result}\label{main_rez}

We intend to determine the Poisson kernel \eqref{classical_Poisson} by listing
some essential properties of the associated integral operator. For that, we use the
normalized $S(U_n\times U_n)$-invariant measure on $S(\mathbb{D})$ for
integration on the Shilov boundary. The principal property of the Poisson kernel is
that the integral operator with this kernel is a morphism of a $SU_{n,n}$-module
of functions on $S(\mathbb{D})$ into a $SU_{n,n}$-module of functions on
$\mathbb{D}$ which takes 1 to 1.

So in the quantum case the required Poisson integral operator is a morphism of the
$U_q\mathfrak{su}_{n,n}$-module $\mathbb{C}[S(\mathbb{D})]_q$ into the
$U_q\mathfrak{su}_{n,n}$-module $\mathscr{D}(\mathbb{D})_q'$ which takes 1 to 1.

Recall that every $u\in\mathscr{D}(\mathbb{D})_q'$ is of the form
$$
u=\sum_{j,k=0}^\infty u_{j,k},\qquad
u_{j,k}\in\mathbb{C}[\operatorname{Mat}_n]_{q,j}
\mathbb{C}[\overline{\operatorname{Mat}}_n]_{q,-k},
$$
and the set
$\left\{z_b^\beta(z_a^\alpha)^*\right\}_{a,b,\alpha,\beta=1,2,\ldots,n}$ is
a basis of the vector space $\mathbb{C}[\operatorname{Mat}_n]_{q,1}
\mathbb{C}[\overline{\operatorname{Mat}}_n]_{q,-1}$. This allows one to
introduce the mixed partial derivatives at zero, the linear functionals
$\left.\dfrac{\partial^2}{\partial
z_b^\beta\partial(z_a^\alpha)^*}\right|_{\mathbf{z}=0}$ such that
$$
u_{1,1}=\sum_{a,b,\alpha,\beta=1}^n\left(\left.\frac{\partial^2u} {\partial
z_b^\beta\partial(z_a^\alpha)^*}\right|_{\mathbf{z}=0}\right)
z_b^\beta\,(z_a^\alpha)^*,\qquad u\in\mathscr{D}(\mathbb{D})_q'.
$$

Now we are in a position to produce a quantum analog of the Hua equations.

\begin{theorem}\label{quantum_Hua_new}
If $u\in\mathscr{D}(\mathbb{D})_q'$ belongs to the image of the Poisson
integral operator on the quantum $n\times n$-matrix ball, then
\begin{equation}\label{q_Hua_new}
\left.\sum_{c=1}^nq^{2c}\frac{\partial^2(\xi u)}{\partial
z_c^\beta\partial(z_c^\alpha)^*}\right|_{\mathbf{z}=0}=0
\end{equation}
for all $\xi\in U_q\mathfrak{sl}_{2n}$, $\alpha,\beta=1,2,\ldots,n$.
\end{theorem}

The equation system \eqref{q_Hua_new} is a $q$-analog of \eqref{classical_Hua_new1}.

\bigskip

\section{Invariant generalized kernels and the associated integral
operators}\label{genker}

List some plausible but less known definitions and results on invariant
integral and integral kernels, see \cite{Vak95}. Consider a Hopf algebra
$A$ and an $A$-module algebra $F$. A linear functional $\nu$ on $F$ is
called an $A$-invariant integral if $\nu$ is a morphism of $A$-modules:
$$\nu(af)=\varepsilon(a)f,\qquad a\in A,\quad f\in F,$$
with $\varepsilon$ being the counit of $A$.

There exists a unique
$U_q\mathfrak{s}(\mathfrak{gl}_n\times\mathfrak{gl}_n)$-invariant integral
$$
\nu:\mathbb{C}[S(\mathbb{D})]_q\to\mathbb{C},\qquad
\nu:\varphi\mapsto\int\limits_{S(\mathbb{D})_q}\varphi d\nu,
$$
which is normalized by $\int\limits_{S(\mathbb{D})_q}1d\nu=1$ (see
\cite[Chapter 3]{Vak01}). It was also demonstrated in \cite{Vak01} that
$\mathbb{C}[S(\mathbb{D})]_q$ is isomorphic to the algebra of regular
functions on the quantum $U_n$ as a
$U_q\mathfrak{s}(\mathfrak{gl}_n\times\mathfrak{gl}_n)$-module $*$-algebra.
This isomorphism can be used to consider the transfer of $\nu$ on the
latter algebra, where it is known to be positive \cite{Woron}. Hence $\nu$
is itself positive: $\int\limits_{S(\mathbb{D})_q}\varphi^*\varphi d\nu>0$
for all non-zero $\varphi\in\mathbb{C}[S(\mathbb{D})]_q$. Our subsequent
results demonstrate how this
$U_q\mathfrak{s}(\mathfrak{gl}_n\times\mathfrak{gl}_n)$-invariant integral
can be used to produce integral operators which are morphisms of
$U_q\mathfrak{sl}_{2n}$-modules.

Consider $A$-module algebras $F_1$, $F_2$. Given a linear functional
$\nu:F_2\to\mathbb{C}$ and $\mathscr{K}\in F_1\otimes F_2$, we associate a
linear integral operator
$$
K:F_2\to F_1,\qquad
K:f\mapsto(\operatorname{id}\otimes\nu)(\mathscr{K}(1\otimes f)).
$$
In this context $\mathscr{K}$ is called the kernel of this integral
operator. Assume that the integral $\nu$ on $F_2$ is invariant and the
bilinear form
\begin{equation*}
f'\times f''\mapsto\nu(f'\;f''),\qquad f',f''\in F_2,
\end{equation*}
is non-degenerate. It is easy to understand that the integral operator with the kernel
$\mathscr{K}$ is a morphism of $A$-modules if and only if this kernel is invariant
\cite{Vak95}. Another statement from \cite{Vak95} that will be used essentially in a
subsequent construction of $A$-invariant kernels is as follows: $A$-invariant kernels
form a subalgebra of $F_1^\mathrm{op}\otimes F_2$, with $F_1^\mathrm{op}$ being the
algebra derived from $F_1$ by replacement of its multiplication by the opposite one.

Additionally to the algebra of kernels
$\operatorname{Pol}(\operatorname{Mat}_n)_q^\mathrm{op}\otimes
\mathbb{C}[S(\mathbb{D})]_q$, we will use the bimodule of generalized
kernels $\mathscr{D}(\mathbb{D}\times S(\mathbb{D}))'_q$ whose elements are
just the formal series
$$
\sum_{i,j}f_{ij}\otimes\varphi_{ij},\qquad
f_{ij}\in\mathbb{C}[\overline{\operatorname{Mat}}_n]_{q,-j}^\mathrm{op}
\mathbb{C}[\operatorname{Mat}_n]_{q,i}^\mathrm{op},\quad
\varphi_{ij}\in\mathbb{C}[S(\mathbb{D})]_q.
$$
This is a bimodule over the algebra
$\operatorname{Pol}(\operatorname{Mat}_n)_q^\mathrm{op}\otimes
\mathbb{C}[S(\mathbb{D})]_q$.

\bigskip

\section{A passage from affine coordinates to homogeneous
coordinates}\label{to-homogeneous}

This Section, just as Sections \ref{prelim} and \ref{genker}, contains some
preliminary material and the known results obtained in \cite{VakSoib90,
SSV4, Vak01}.

\bigskip

Turn to the quantum group $SL_{N}$. We follow the general idea of Drinfeld
\cite{Drinf1} in considering the Hopf algebra $\mathbb{C}[SL_{N}]_q$ of
matrix elements of finite dimensional weight
$U_q\mathfrak{sl}_{N}$-modules. It is custom to call it the `algebra of
regular functions on the quantum group $SL_{N}$'. The linear maps in
$(U_q\mathfrak{sl}_{N})^*$ adjoint to the operators of left multiplication
by the elements of $U_q\mathfrak{sl}_{N}$ equip $\mathbb{C}[SL_{N}]_q$ with
a structure of $U_q\mathfrak{sl}_{N}$-module algebra by duality.

Recall that $\mathbb{C}[SL_{N}]_q$ can be defined by the generators
$t_{ij}$, $i,j=1,...,N$, (the matrix elements of the vector representation
in a weight basis) and the relations
\begin{align*}
&t_{ij'}t_{ij''}=qt_{ij''}t_{ij'},\qquad &j'<j'',
\\ &t_{i'j}t_{i''j}=qt_{i''j}t_{i'j},\qquad &i'<i'',
\\ &t_{ij}t_{i'j'}=t_{i'j'}t_{ij},\qquad &i<i'\;\&\;j>j',
\\ &t_{ij}t_{i'j'}=t_{i'j'}t_{ij}+(q-q^{-1})t_{ij'}t_{i'j},\qquad
& i<i'\;\&\; j<j',
\end{align*}
which are tantamount to \eqref{z}, together with one more relation
\begin{equation*}\label{det_1}
\det\nolimits_q\mathbf{t}=1,
\end{equation*}
where $\det\nolimits_q\mathbf{t}$ is a $q$-determinant of the matrix
$\mathbf{t}=(t_{ij})_{i,j=1,...,N}$:
\begin{equation*}
\det\nolimits_q\mathbf{t}=\sum\limits_{s\in S_{N}}(-q)^{l(s)}t_{1 s(1)}
t_{2s(2)}\ldots t_{Ns(N)},
\end{equation*}
with $l(s)=\mathrm{card}\{(i,j)|i<j\;\&\;s(i)>s(j)\}$.

It is well known that $\det\nolimits_q\mathbf{t}$ commutes with all
$t_{ij}$. Thus $\mathbb{C}[SL_{N}]_q$ appears to be a quotient algebra of
$\mathbb{C}[\mathrm{Mat}_{N}]_q$ by the two-sided ideal generated by
$\mathrm{det}_q\mathbf{t}-1$.

Note that $\mathbb{C}[SL_{N}]_q$ is a domain.

\bigskip

In the classical case $q=1$ the matrix ball admits a natural embedding into
the Grassmannian $\mathrm{Gr}_{n,2n}$
$$
(\text{a contraction\ }A\in\operatorname{End}\mathbb{C}^n)\mapsto
(\text{the linear span of\ }(v,Av),\quad v\in\mathbb{C}^n),
$$
with the latter pair being an element of
$\mathbb{C}^n\oplus\mathbb{C}^n\simeq\mathbb{C}^{2n}$. We are going to
describe a $q$-analog for this embedding.

Let
$$
I=\{i_1,i_2,\ldots,i_k\}\subset\{1,2,\ldots,2n\},\qquad{i_1<i_2<\ldots<i_k};
$$
$$
J=\{j_1,j_2,\ldots,j_k\}\subset\{1,2,\ldots,2n\},\qquad j_1<j_2<\ldots<j_k.
$$
The elements
$$
t_{IJ}^{\wedge k}=\sum\limits_{s\in S_k}(-q)^{l(s)}t_{i_{s(1)}j_1}
t_{i_{s(2)}j_2}\ldots t_{i_{s(k)}j_k}
$$
of $\mathbb{C}[SL_{2n}]_q$ are called quantum minors, and it is easy to
check that
\begin{equation*}\label{q-minor_new}
t_{IJ}^{\wedge k}=\sum\limits_{s\in S_k}(-q)^{l(s)}t_{i_1j_{s(1)}}
t_{i_2j_{s(2)}}\ldots t_{i_k j_{s(k)}}.
\end{equation*}

Consider the smallest unital subalgebra
$\mathbb{C}[X]_q\subset\mathbb{C}[SL_{2n}]_q$ that contains the quantum
minors
$$
t_{\{1,2,\ldots,n\}J}^{\wedge n},\;\;\;t_{\{n+1,n+2,\ldots,2n\}J}^{\wedge
n},\qquad J=\{j_1,j_2,\ldots,j_n\}\subset\{1,2,\ldots,2n\}.
$$
It is a $U_q\mathfrak{sl}_{2n}$-module subalgebra which substitutes the
classical coordinate ring of the Grassmannian.

The following results are easy modifications of those of \cite{SSV4}.

\begin{proposition}
There exists a unique antilinear involution $*$ in $\mathbb{C}[X]_q$ such
that $(\mathbb{C}[X]_q,*)$ is a $U_q\mathfrak{su}_{n,n}$-module algebra and
$$
\left(t_{\{1,2,\ldots,n\}\{n+1,n+2,\ldots,2n\}}^{\wedge n}\right)^*=
(-q)^{n^2}t_{\{n+1,n+2,\ldots,2n\}\{1,2,\ldots,n\}}^{\wedge n}.
$$
\end{proposition}

\begin{lemma}[Lemma 11.3 of \cite{SSV4}]\label{invfl}
Given $J\subset\{1,2,\ldots,2n\}$ with $\mathrm{card}(J)=n$,
\hbox{$J^c=\{1,2,\ldots,2n\}\setminus J$}, $l(J,J^c)=\mathrm{card}\{(j',j'')\in
J\times J^c|\:j'>j''\}$. Then
\begin{equation}\label{invf}
\left(t^{\wedge
n}_{\{1,2,\ldots,n\}J}\right)^*=(-1)^{\mathrm{card}(\{1,2,\ldots,n\}\cap
J)}(-q)^{l(J,J^c)}t^{\wedge n}_{\{n+1,n+2,\ldots,2n\}J^c}.
\end{equation}
\end{lemma}

Impose the abbreviated notation
$$
t=t_{\{1,2,\ldots,n\}\{n+1,n+2,\ldots,2n\}}^{\wedge n},\qquad x=tt^*.
$$
Note that $t$, $t^*$ and $x$ quasi-commute with all the generators $t_{ij}$
of $\mathbb{C}[SL_{2n}]_q$. Then the localization $\mathbb{C}[X]_{q,x}$ of
the algebra $\mathbb{C}[X]_q$ with respect to the multiplicative set
$x^{\mathbb{Z}_+}$ is well-defined. The structure of
$U_q\mathfrak{su}_{n,n}$-module algebra is uniquely extendable from
$\mathbb{C}[X]_q$ onto $\mathbb{C}[X]_{q,x}$ \cite{SSV4}.

\begin{proposition}[cf. Proposition 3.2 of \cite{SSV4}]
There exists a unique embedding of $U_q\mathfrak{su}_{n,n}$-module
$*$-algebras
$\mathcal{I}:\operatorname{Pol}(\operatorname{Mat}_n)_q\hookrightarrow
\mathbb{C}[X]_{q,x}$ such that $\mathcal{I}z_a^\alpha=t^{-1}t^{\wedge
n}_{\{1,2,\ldots,n\}J_{a\alpha}}$, with
$J_{a\alpha}=\{a\}\cup\{n+1,n+2,\ldots,2n\}\backslash\{2n+1-\alpha\}$.
\end{proposition}

\begin{corollary}
$\mathcal{I}y=x^{-1}$, with
\begin{equation}\label{y_mat}
y=1+\sum\limits_{k=1}^m (-1)^k
\sum\limits_{\{J'|\:\operatorname{card}(J')=k\}}
\sum\limits_{\{J''|\:\operatorname{card}(J'')=k\}}z_{\quad J''}^{\wedge
k\;J'}\left(z_{\quad J''}^{\wedge k\;J'}\right)^*.
\end{equation}
\end{corollary}
A formal passage to a limit as $q\to 1$ leads to the relation
$y=\det(1-\mathbf{z}\mathbf{z}^*)$.

\begin{proposition}[see Chapter 11 of \cite{SSV4}]
Let $1\le\alpha_1<\alpha_2<\ldots<\alpha_k\le n$, $1\le
a_1<a_2<\ldots<a_k\le n$,
$J=\{n+1,n+2,\ldots,2n\}\setminus\{n+\alpha_1,n+\alpha_2,\ldots,n+\alpha_k\}
\cup\{a_1,a_2,\ldots,a_k\}$. Then
\begin{equation}\label{icmat}
\mathcal{I}z_{\phantom{\wedge k}\,\{a_1,a_2,\ldots,a_k\}}^{\wedge
k\{n+1-\alpha_k,n+1-\alpha_{k-1},\ldots,
n+1-\alpha_1\}}=t^{-1}t_{\{1,2,\ldots,n\}J}^{\wedge n}.
\end{equation}
\end{proposition}

It is custom to identify the generators $z_a^\alpha$,
$a,\alpha=1,2,\ldots,n$, with their images under $\mathcal{I}$.

Consider the subalgebra of $\mathbb{C}[X]_{q,x}$ generated by
$\operatorname{Pol}(\operatorname{Mat}_n)_q
\overset{\mathcal{I}}{\hookrightarrow}\mathbb{C}[X]_{q,x}$, together with
$t^{\pm 1}$, $t^{*\pm 1}$.

The elements of this subalgebra admit a unique decomposition of the form
\begin{equation}\label{decomp}
\sum_{(i,j)\notin(-\mathbb{N})\times(-\mathbb{N})}t^it^{*j}f_{ij},\qquad
f_{ij}\in\operatorname{Pol}(\operatorname{Mat}_n)_q
\end{equation}
(the choice of the set of pairs $(i,j)$ is due to the fact that
$t^{-1}t^{*-1}\in\operatorname{Pol}(\operatorname{Mat}_n)_q$).

\bigskip

Equip $\mathbb{C}^n\oplus\mathbb{C}^n$ with the sesquilinear form
$(\cdot,\cdot)_1-(\cdot,\cdot)_2$. In the classical case $q=1$ the Shilov
boundary of the matrix ball is just the group $U_n$. The graphs of unitary
operators form the isotropic Grassmannian (its points are just the
subspaces on which the above scalar product on $\mathbb{C}^{2n}$ vanishes).
We are going to describe a $q$-analog of it, with the isotropic
Grassmannian being replaced by its homogeneous coordinate ring.

Consider an extension $\mathbb{C}[\Xi]_q$ of the algebra
$\mathbb{C}[S(\mathbb{D})]_q$ in the class of
$U_q\mathfrak{su}_{n,n}$-module $*$-algebras. This extension is produced
via adding a generator $t$ and the relations
\begin{equation}\label{tz}
tt^*=t^*t;\quad tz_a^\alpha=q^{-1}z_a^\alpha t;\quad t^*z_a^\alpha=q^{-1}z_a^\alpha
t^*,\quad a,\alpha=1,2,\ldots,n.
\end{equation}
The $U_q\mathfrak{su}_{n,n}$-action is extended onto $\mathbb{C}[\Xi]_q$ as
follows:
$$E_jt=F_jt=(K_j^{\pm 1}-1)t=0,\qquad j\ne n,$$
$$F_nt=(K_n^{\pm 1}-1)t=0,\qquad E_nt=q^{-1/2}tz_n^n.$$

Let $\xi=tt^*$. One can introduce a localization $\mathbb{C}[\Xi]_{q,\xi}$
of the algebra $\mathbb{C}[\Xi]_q$ with respect to the multiplicative set
$\xi^{\mathbb{Z}_+}$. The involution $*$ and the structure of
$U_q\mathfrak{su}_{n,n}$-module algebra are uniquely extendable from
$\mathbb{C}[\Xi]_q$ onto $\mathbb{C}[\Xi]_{q,\xi}$.

The algebra $\mathbb{C}[\Xi]_q$ is equipped with a
$U_q\mathfrak{sl}_{2n}$-invariant bigrading
$$
\deg z_a^\alpha=\deg(z_a^\alpha)^*=(0,0),\qquad\deg t=(1,0),\qquad\deg
t^*=(0,1),
$$
which extends onto the localization $\mathbb{C}[\Xi]_{q,\xi}$. It was
demonstrated in \cite{Vak01} that the homogeneous component
$\mathbb{C}[\Xi]_{q,\xi}^{(-n,-n)}$ carries a non-zero
$U_q\mathfrak{sl}_{2n}$-invariant integral $\eta$ such that
\begin{equation*}
\eta(t^{*-n}ft^{-n})=\int\limits_{S(\mathbb{D})_q}fd\nu.
\end{equation*}
Since $t$, $t^*$ normalize every
$U_q\mathfrak{s}(\mathfrak{u}_n\times\mathfrak{u}_n)$-isotypical component
of $\mathbb{C}[S(\mathbb{D})]_q$ and the
$U_q\mathfrak{s}(\mathfrak{u}_n\times\mathfrak{u}_n)$-action in
$\mathbb{C}[S(\mathbb{D})]_q$ is multiplicity free \cite[Proposition
4]{my_deg}, it follows from the definition of the integral over the Shilov
boundary that
\begin{equation}\label{qcic}
\eta(t^{*-n}t^{-n}f)=\eta(ft^{*-n}t^{-n})=
\int\limits_{S(\mathbb{D})_q}fd\nu.
\end{equation}

\bigskip

\section{The Poisson kernel}\label{Poisson-Shilov}

In this Section we construct a morphism of the
$U_q\mathfrak{su}_{n,n}$-module $\mathbb{C}[S(\mathbb{D})]_q$ into the
$U_q\mathfrak{su}_{n,n}$-module $\mathscr{D}(\mathbb{D})_q'$ which takes 1
to 1. The introduced morphism will be a Poisson integral operator, as
mentioned in Section \ref{main_rez}.

We will need a quantum analog $\mathbb{C}[\operatorname{Mat}_{n,2n}]_q$ for
the polynomial algebra on the space of rectangular $n\times 2n$-matrices.
The algebra $\mathbb{C}[\operatorname{Mat}_{n,2n}]_q$ is determined by the
set of generators $\{t_{ij}\}_{i=1,2,\ldots,n;\;j=1,2,\ldots,2n}$ and the
similar relations as in the case of square matrices, see \eqref{z}. It is
known that this algebra is a domain.

The structure of $U_q\mathfrak{sl}_{2n}$-module algebra is given by
$$
K_kt_{ij}=\left\{
\begin{array}{rl}
qt_{ij}, & \quad j=k,
\\ q^{-1}t_{ij}, & \quad j=k+1,
\\ t_{ij}, & \quad j\notin\{k,k+1\},
\end{array}\right.
$$
$$
E_kt_{ij}=\left\{
\begin{array}{rl}
q^{-1/2}t_{i,j-1}, & j=k+1,
\\ 0, & j\ne k+1,
\end{array}\right.\qquad
F_kt_{ij}=\left\{
\begin{array}{rl}
q^{1/2}t_{i,j+1}, & j=k,
\\ 0, & j\ne k,
\end{array}\right.\
$$
with ${k=1,2,\ldots,2n-1}$.

\begin{remark}\label{embs}
Consider the $U_q\mathfrak{sl}_{2n}$-module subalgebras in
$\mathbb{C}[\operatorname{Mat}_{n,2n}]_q$ and in $\mathbb{C}[\Xi]_q$ generated by
$t^{\wedge n}_{\{1,2,\ldots,n\}\{n+1,n+2,\ldots,2n\}}$ and $t$, respectively. The map
$t^{\wedge n}_{\{1,2,\ldots,n\}\{n+1,n+2,\ldots,2n\}}\mapsto t$ extends up to an
isomorphism $i$ of these $U_q\mathfrak{sl}_{2n}$-module subalgebras. In a similar way,
one can introduce $U_q\mathfrak{sl}_{2n}$-module subalgebras in
$\mathbb{C}[\operatorname{Mat}_{n,2n}]_q$ and in $\mathbb{C}[\Xi]_q$, generated by
$t^{\wedge n}_{\{1,2,\ldots,n\}\{1,2,\ldots,n\}}$ and $t^*$, respectively. The map
$t^{\wedge n}_{\{1,2,\ldots,n\}\{1,2,\ldots,n\}}\mapsto(-q)^{n^2}t^*$ extends up to an
isomorphism $i'$ of these $U_q\mathfrak{sl}_{2n}$-module subalgebras. The above
isomorphisms can be used together to embed minors $t^{\wedge n}_{\{1,2,\ldots,n\}J}$ and
$t^{\wedge n}_{\{n+1,n+2,\ldots,2n\}J}$ into $\mathbb{C}[\Xi]_q$.
\end{remark}

Consider the map $m:\mathbb{C}[\operatorname{Mat}_{n,2n}]_q\otimes
\mathbb{C}[\operatorname{Mat}_{n,2n}]_q\to \mathbb{C}[\operatorname{Mat}_{2n}]_q$ defined
as follows. The tensor multipliers of the domain are embedded into
$\mathbb{C}[\operatorname{Mat}_{2n}]_q$ as subalgebras generated by the entries of,
respectively, the upper $n$ and the lower $n$ rows of the matrix $\mathbf{t}=(t_{ij})$,
$i=1,\ldots,n$, $j=1,\ldots,2n$, and the map $m$ is just the multiplication in the
algebra $\mathbb{C}[\operatorname{Mat}_{2n}]_q$.

\begin{lemma}\label{m_iso}
The map $m$ is an isomorphism of $U_q\mathfrak{sl}_{2n}$-modules.
\end{lemma}

{\bf Proof.} The map $m$ takes the monomial basis
$$
t_{11}^{j_{11}}\ldots t_{1,2n}^{j_{1,2n}}t_{21}^{j_{21}}\ldots
t_{2,2n}^{j_{2,2n}}\ldots t_{n1}^{j_{n1}}\ldots t_{n,2n}^{j_{n,2n}}\otimes
t_{n+1,1}^{j_{n+1,1}}\ldots
t_{n+1,2n}^{j_{n+1,2n}}t_{n+2,1}^{j_{n+2,1}}\ldots
t_{n+2,2n}^{j_{n+2,2n}}\ldots t_{2n,1}^{j_{2n,1}}\ldots
t_{2n,2n}^{j_{2n,2n}}
$$
of the algebra $\mathbb{C}[\operatorname{Mat}_{n,2n}]_q\otimes
\mathbb{C}[\operatorname{Mat}_{n,2n}]_q$ to the monomial basis
$$
t_{11}^{j_{11}}\ldots t_{1,2n}^{j_{1,2n}}t_{21}^{j_{21}}\ldots
t_{2,2n}^{j_{2,2n}}\ldots t_{n1}^{j_{n1}}\ldots t_{n,2n}^{j_{n,2n}}
t_{n+1,1}^{j_{n+1,1}}\ldots t_{n+1,2n}^{j_{n+1,2n}}\cdot
t_{n+2,1}^{j_{n+2,1}}\ldots t_{n+2,2n}^{j_{n+2,2n}}\ldots
t_{2n,1}^{j_{2n,1}}\ldots t_{2n,2n}^{j_{2n,2n}}
$$
of the algebra $\mathbb{C}[\operatorname{Mat}_{2n}]_q$,
$j_{ik}\in\mathbb{Z}_+$, hence is a bijective map. It is a morphism of
$U_q\mathfrak{sl}_{2n}$-modules since
$\mathbb{C}[\operatorname{Mat}_{2n}]_q$ is a $U_q\mathfrak{sl}_{2n}$-module
algebra. \hfill $\square$

\medskip

It is worthwhile to note that the definition of the
$U_q\mathfrak{sl}_{2n}$-module algebra $\mathbb{C}[X]_q$ allows a
replacement of $\mathbb{C}[SL_{2n}]_q$ by
$\mathbb{C}[\operatorname{Mat}_{n,2n}]_q$.

Consider the elements of the $U_q\mathfrak{sl}_{2n}$-module
$\mathbb{C}[\operatorname{Mat}_{n,2n}]_q\otimes
\mathbb{C}[\operatorname{Mat}_{n,2n}]_q$ given by
$$
\mathscr{L}=
\sum\limits_{J\subset\{1,2,\ldots,2n\}\&\operatorname{card}(J)=n}
(-q)^{l(J,J^c)}t^{\wedge n}_{\{1,2,\ldots,n\}J}\otimes t^{\wedge
n}_{\{n+1,n+2,\ldots,2n\}J^c},
$$
$$
\overline{\mathscr{L}}=
\sum\limits_{J\subset\{1,2,\ldots,2n\}\&\mathrm{card}(J)=n}
(-q)^{-l(J,J^c)}t^{\wedge n}_{\{n+1,n+2,\ldots,2n\}J^c}\otimes t^{\wedge
n}_{\{1,2,\ldots,n\}J}.
$$
Here $J^c$ is the complement to $J$ and
$l(I,J)=\operatorname{card}\{(i,j)\in I\times J|\:i>j\}$.

\begin{proposition}
$\mathscr{L}$ and $\overline{\mathscr{L}}$ are
$U_q\mathfrak{sl}_{2n}$-invariants.
\end{proposition}

{\bf Proof} is expounded for $\mathscr{L}$. In the case of
$\overline{\mathscr{L}}$, similar arguments are applicable.

Recall a $q$-analog for the Laplace formula of splitting the quantum
determinant of the $2n\times 2n$-matrix $\mathbf{t}=(t_{ij})$ with respect
to the upper $n$ lines:
\begin{multline*}
\det\nolimits_q\mathbf{t}=
\sum\limits_{J\subset\{1,2,\ldots,2n\}\&\operatorname{card}(J)=n}
(-q)^{l(J,J^c)}t^{\wedge n}_{\{1,2,\ldots,n\}J}t^{\wedge
n}_{\{n+1,n+2,\ldots,2n\}J^c}=
\\ =\sum\limits_{J\subset\{1,2,\ldots,2n\}\&\mathrm{card}(J)=n}
(-q)^{-l(J,J^c)}t^{\wedge n}_{\{n+1,n+2,\ldots,2n\}J^c}t^{\wedge
n}_{\{1,2,\ldots,n\}J}.
\end{multline*}
Our claim follows from Lemma \ref{m_iso}, the relation
$m\mathscr{L}=\det_q\mathbf{t}$ and $U_q\mathfrak{sl}_{2n}$-invariance of
the quantum determinant. \hfill $\square$

\medskip

Note that, in view of Remark \ref{embs}, one has
$\mathscr{L},\overline{\mathscr{L}}\in\mathbb{C}[X]_q\otimes
\mathbb{C}[\Xi]_q$.

Introduce, firstly, a $U_q\mathfrak{sl}_{2n}$-module of kernels
$\mathscr{D}(\mathbb{D}\times\Xi)'_q$, whose elements are formal series
with coefficients from $\mathbb{C}[\overline{\operatorname{Mat}}_n]_{q,-j}
\mathbb{C}[\operatorname{Mat}_n]_{q,i}\otimes\mathbb{C}[\Xi]_{q,\xi}$, and,
secondly, a $U_q\mathfrak{sl}_{2n}$-module of kernels
$\mathscr{D}(X\otimes\Xi)'_q$, whose elements are finite sums of the form
\begin{equation*}
\sum_{(i,j)\notin(-\mathbb{N})\times(-\mathbb{N})}(t^it^{*j}\otimes
1)f_{ij},\qquad f_{ij}\in\mathscr{D}(\mathbb{D}\times\Xi)'_q
\end{equation*}
(cf. \eqref{decomp}). Of course, $\mathscr{D}(\mathbb{D}\times\Xi)'_q$ is a
$\operatorname{Pol}(\operatorname{Mat}_n)_q^\mathrm{op}\otimes\
\mathbb{C}[\Xi]_{q,\xi}$-bimodule, and $\mathscr{D}(X\times\Xi)'_q$ is a
$\mathbb{C}[X]_{q,x}^\mathrm{op}\otimes\mathbb{C}[\Xi]_{q,\xi}$-bimodule.

\bigskip \bigskip \bigskip

The kernel
$\mathscr{L}\in\mathbb{C}[X]_q^\mathrm{op}\otimes\mathbb{C}[\Xi]_q$ can be
written in the form
$$
\mathscr{L}=(-q)^{-n^2}\cdot
\left(1+\sum\limits_{J\ne\{1,2,\ldots,n\}}(-q)^{l(J,J^c)}
t_{\{1,2,\ldots,n\}J}^{\wedge n}t^{-1}\otimes t^{\wedge
n}_{\{n+1,n+2,\ldots,2n\}J^c}t^{*-1}\right)(t\otimes t^*).
$$
Note that $t_{\{1,2,\ldots,n\}J}^{\wedge n}t^{-1},t^{\wedge
n}_{\{n+1,n+2,\ldots,2n\}J^c}t^{*-1}\in\mathbb{C}[\operatorname{Mat}_n]_q$,
see \eqref{icmat}. This allows one to write down explicitly such element
$\mathscr{L}^{-n}$ of the space of generalized kernels
$\mathscr{D}(X\times\Xi)'_q$ that
$\mathscr{L}^n\cdot\mathscr{L}^{-n}=\mathscr{L}^{-n}\cdot\mathscr{L}^n=1$,
where $\cdot$ stands for the (left and right) actions of $\mathscr{L}^n$ on
the element $\mathscr{L}^{-n}$ of the bimodule
$\mathscr{D}(X\times\Xi)'_q$.

A similar construction produces also a $U_q\mathfrak{sl}_{2n}$-invariant generalized
kernel $\overline{\mathscr{L}}^{\,-n}$.

Note that $\mathscr{L}^{-n}=\sum\limits_ix_i$ is a formal series with
$x_i\in\mathbb{C}[\operatorname{Mat}_n]_{q,i}\otimes\mathbb{C}[\Xi]_{q,\xi}$,
$i\in\mathbb{Z}_+$, and the terms of the formal series
$\overline{\mathscr{L}}^{\,-n}=\sum\limits_jy_j$ are such that
$y_j\in\mathbb{C}[\overline{\operatorname{Mat}}_n]_{q,-j}
\otimes\mathbb{C}[\Xi]_{q,\xi}$. This allows one to define the `product'
$\overline{\mathscr{L}}^{\,-n}\mathscr{L}^{-n}$ as a double series
$\sum\limits_{i,j}y_jx_i$ which is thus an element of the module of
generalized kernels $\mathscr{D}(X\otimes\Xi)'_q$. Clearly, one has
$\overline{\mathscr{L}}^n\cdot
\left(\overline{\mathscr{L}}^{\,-n}\mathscr{L}^{-n}\right)
\cdot\mathscr{L}^n=1$ in $\mathscr{D}(X\times\Xi)'_q$, and this property
determines uniquely the generalized kernel
$\overline{\mathscr{L}}^{\,-n}\mathscr{L}^{-n}$. Furthermore, the above
uniqueness allows one to verify invariance of the generalized kernel
$\overline{\mathscr{L}}^{\,-n}\mathscr{L}^{-n}$. Of course the argument
should apply the invariance of $\mathscr{L}$, $\overline{\mathscr{L}}$.

\bigskip

Consider the Poisson kernel $P\in\mathscr{D}(\mathbb{D}\times
S(\mathbb{D}))'_q$ for the matrix ball by demanding the following
properties:

i) up to a constant multiplier, the Poisson kernel is just $(1\otimes
tt^*)^{n}\overline{\mathscr{L}}^{\,-n}\mathscr{L}^{-n}$, that is
\begin{equation*}\label{Poisson}
P=\mathrm{const}(q,n)(1\otimes
tt^*)^{n}\overline{\mathscr{L}}^{\,-n}\mathscr{L}^{-n};
\end{equation*}

ii) the integral operator with kernel $P$ takes
$1\in\mathbb{C}[S(\mathbb{D})]_q$ to $1\in\mathscr{D}(\mathbb{D})'_q$.

\medskip

\begin{lemma}\label{iio}
The integral operator
$\mathbb{C}[S(\mathbb{D})]_q\to\mathscr{D}(\mathbb{D})'_q$ with
kernel $(1\otimes
tt^*)^{n}\overline{\mathscr{L}}^{\,-n}\mathscr{L}^{-n}$ is a morphism
of $U_q\mathfrak{sl}_{2n}$-modules.
\end{lemma}

{\bf Proof.} It follows from the existence of an invariant integral
$\eta:\mathbb{C}[\Xi]_q^{(-n,-n)}\to\mathbb{C}$ and the invariance of
$\overline{\mathscr{L}}^{\,-n}\mathscr{L}^{-n}$ that the integral operator
$\mathbb{C}[S(\mathbb{D})]_q\to\mathscr{D}(\mathbb{D})'_q$ with kernel
$\overline{\mathscr{L}}^{\,-n}\mathscr{L}^{-n}$ is a morphism of
$U_q\mathfrak{sl}_{2n}$-modules. Now it remains to move the multiplier
$(1\otimes t^{*-n}t^{-n})$ from
$\overline{\mathscr{L}}^{\,-n}\mathscr{L}^{-n}$ to the left and to apply
\eqref{qcic}.\hfill $\square$

\begin{remark}
Since the integral operator
$\mathbb{C}[S(\mathbb{D})]_q\to\mathscr{D}(\mathbb{D})'_q$ with
kernel $(1\otimes
tt^*)^{n}\overline{\mathscr{L}}^{\,-n}\mathscr{L}^{-n}$ is a morphism
of $U_q\mathfrak{sl}_{2n}$-modules, the image of $1 \in
\mathbb{C}[S(\mathbb {D})]_q$ is a $U_q \mathfrak{sl}_{2n}$-invariant
element of $\mathscr{D}(\mathbb{D})'_q$, that is just a constant.
This proves the existence of the Poisson kernel $P$.
\end{remark}

\medskip

\begin{example}
Let us illustrate the Poisson kernel $P$ determined above in the simplest case $n=1$. The
considered invariant kernels were obtained in the early paper \cite{ftqdik}. The elements
$$
\mathscr{L}=t_{11}\otimes t_{22}-qt_{12}\otimes t_{21},\qquad
\overline{\mathscr{L}}=-q^{-1}t_{21}\otimes t_{12}+t_{22}\otimes t_{11}
$$
of the algebra
$\mathbb{C}[X]_{q,x}^\mathrm{op}\otimes\mathbb{C}[\Xi]_{q,\xi}$ are
$U_q\mathfrak{sl}_2$-invariant kernels. We present below an easy
computation, which uses, for the sake of brevity, the notation
\begin{equation*}
z=qt_{12}^{-1}t_{11},\qquad z^*=t_{22}t_{21}^{-1},
\end{equation*}
instead of mentioning explicitly the embeddings $\mathcal{I}$,
$\mathscr{I}$. One has:
$$
\mathscr{L}=(1-z\otimes z^*)(-qt_{12}\otimes t_{21}),\qquad
\overline{\mathscr{L}}=(-q^{-1}t_{21}\otimes t_{12})(1-q^2z^*\otimes z),
$$
and, as $t=t_{12}$, $t^*=-qt_{21}$,
$$
\mathscr{L}=(1-z\otimes z^*)(t\otimes t^*),\qquad
\overline{\mathscr{L}}=(q^{-2}t^*\otimes t)(1-q^2z^*\otimes z).
$$
Hence,
$$
\overline{\mathscr{L}}^{\,-1}\mathscr{L}^{-1}=q^2(1\otimes
t^{-1}t^{*-1})(1-z^*\otimes z)^{-1}((1-z^*z)\otimes 1)(1-z\otimes
z^*)^{-1}.
$$
Omit $\otimes$ and replace in the second tensor multiplier $z$ by $\zeta$
and $z^*$ by $\zeta^*$ (which is standard in function theory) to obtain
\begin{equation*}
P=\mathrm{const}(q)(1-z^*\zeta)^{-1}(1-z^*z)(1-z\zeta^*)^{-1}.
\end{equation*}
What remains now is to find $\mathrm{const}(q)$, or, to be more precise, to
prove that $\mathrm{const}(q)=1$. In fact, the integral operator with
kernel $(1-z^*\zeta)^{-1}(1-z^*z)(1-z\zeta^*)^{-1}$ takes $1$ to $1$. This
is because $\zeta^*=\zeta^{-1}$, and integration of the product of the
series in $\zeta$ produces the constant term: $\sum\limits_{k=0}^\infty
z^k(1-zz^*)z^{*k}=1$.
\end{example}

\medskip

Note that the Poisson kernel is a formal series
$P=\sum\limits_{j,k=0}^\infty p_{jk}$, with
$p_{jk}\in\mathbb{C}[\overline{\operatorname{Mat}}_n]_{q,-k}
\mathbb{C}[\operatorname{Mat}_n]_{q,j}\otimes\mathbb{C}[S(\mathbb{D})]_q$.
In the sequel we will omit $\otimes$ and replace in the second tensor
multiplier $z$ by $\zeta$ and $z^*$ by $\zeta^*$ (which is standard in
function theory).

\begin{lemma}\label{P_q}
The following relation is valid:
\begin{equation}\label{p_1_1}
p_{11}=\mathrm{const}(q,n)
\sum_{a,b,\alpha,\beta=1}^n\left(\frac{1-q^{-2n}}{1-q^{-2}}
q^{2(2n-a-\alpha)}\zeta_a^\alpha(\zeta_b^\beta)^*
-\delta_{ab}\delta^{\alpha\beta}\right)(z_a^\alpha)^*z_b^\beta,
\end{equation}
with $\mathrm{const}(q,n)\ne 0$.
\end{lemma}

{\bf Proof.} In the algebra of kernels one has
\begin{equation}\label{L_up_to}
\mathscr{L}=\left(1-\sum_{a,\alpha=1}^nz_a^\alpha(\zeta_a^\alpha)^*+
\ldots\right)t\tau^*,
\end{equation}
\begin{equation}\label{Lbar_up_to}
\overline{\mathscr{L}}=q^{-2n^2}t^*\tau\left(1-q^2\sum_{a,\alpha=1}^n
q^{2(2n-a-\alpha)}(z_a^\alpha)^*\zeta_a^\alpha+\ldots\right),
\end{equation}
which is easily deducible from \eqref{icmat}, \eqref{invf}, and the fact
that $\mathcal{I}$ is a homomorphism of $*$-algebras. Also, it is very well
visible from \eqref{y_mat} that
\begin{equation}\label{y_up_to}
y=(tt^*)^{-1}=1-\sum_{a,\alpha=1}^n(z_a^\alpha)^*z_a^\alpha+\ldots.
\end{equation}
Here three dots replace the terms whose degree is above two, and the
following abbreviated notation is implicit:
$$
t=t\otimes 1,\qquad\tau=1\otimes\tau,\qquad z_a^\alpha=z_a^\alpha\otimes
1,\qquad\zeta_a^\alpha=1\otimes\zeta_a^\alpha,
$$
$$
t^*=t^*\otimes 1,\qquad\tau^*=1\otimes\tau^*,\qquad(z_a^\alpha)^*=
(z_a^\alpha)^*\otimes
1,\qquad(\zeta_a^\alpha)^*=1\otimes(\zeta_a^\alpha)^*.
$$
Apply \eqref{L_up_to} -- \eqref{y_up_to}, together with the commutation
relations (see also \eqref{tz})
$$
t^*\tau\left(\sum_{a,\alpha=1}^nq^{2(2n-a-\alpha)}
(z_a^\alpha)^*\zeta_a^\alpha\right)=q^{-2}\left(\sum_{a,\alpha=1}^n
q^{2(2n-a-\alpha)}(z_a^\alpha)^*\zeta_a^\alpha\right)t^*\tau
$$
$$
t\tau^*\left(\sum_{a,\alpha=1}^nz_a^\alpha(\zeta_a^\alpha)^*\right)=
q^2\left(\sum_{a,\alpha=1}^nz_a^\alpha(\zeta_a^\alpha)^*\right)t\tau^*,
$$
to obtain:
\begin{eqnarray*}
\overline{\mathscr{L}}^{\,-n}\mathscr{L}^{-n} &=& \prod_{j=1}^n
\left(1-q^{2j}\sum_{a,\alpha=1}^nq^{2(2n-a-\alpha)}
(z_a^\alpha)^*\zeta_a^\alpha+\ldots\right)^{-1}\cdot\nonumber
\\ & & \cdot q^{2n^3}(\tau^*\tau)^{-n}(tt^*)^{-n}\prod_{j=0}^{n-1}
\left(1-q^{2j}\sum_{a,\alpha=1}^nz_a^\alpha(\zeta_a^\alpha)^*+\ldots
\right)^{-1}.
\end{eqnarray*}
Hence
\begin{multline}\label{P_up_to}
P=\mathrm{const}(q,n)\prod_{j=0}^{n-1}\left(1-q^{-2j}
\sum_{a,\alpha=1}^nq^{2(2n-a-\alpha)}(z_a^\alpha)^*\zeta_a^\alpha+
\ldots\right)^{-1}\cdot
\\  \cdot\left(1-\sum_{a,\alpha=1}^n(z_a^\alpha)^*z_a^\alpha+\ldots
\right)^n\cdot\prod_{j=0}^{n-1}\left(1-q^{2j}
\sum_{a,\alpha=1}^nz_a^\alpha(\zeta_a^\alpha)^*+\ldots\right)^{-1},
\end{multline}
where three dots inside every parentheses denotes the contribution of the
terms whose degree is above two, and the multipliers in the products are
written in order of decrease of the index $j$ from left to right.

Now \eqref{p_1_1} follows from \eqref{P_up_to}. \hfill $\square$

\begin{remark}\label{classical_p_1_1}
A formal passage to a limit as $q\to 1$ in \eqref{p_1_1} leads to
\begin{equation*}
p_{11}=\mathrm{const}(n)\sum_{a,b,\alpha,\beta=1}^n \left(n\,\zeta_a^\alpha
\overline{\zeta_b^\beta}-
\delta_{ab}\delta^{\alpha\beta}\right)\overline{z_a^\alpha}z_b^\beta.
\end{equation*}
This relation is well known (with $\mathrm{const}(n)=n$), see, for example,
\cite[p. 597]{JohnsonKoranyi}, and is a consequence of
\eqref{classical_Poisson}.
\end{remark}
\begin{remark}
Using the definition of $P$, \eqref{P_up_to} and the definition of
the integral over the Shilov boundary, one can compute
$\mathrm{const}(q,n)$ explicitly. But we do not need this value on
the way to Hua equations.
\end{remark}

\bigskip

\section{Deducing the Hua equations}\label{Hua_system}

Now we are about to produce a quantum analog of \eqref{Hua_zero_classical1}.

It follows from Lemma \ref{P_q} and the definition of multiplication in the
algebra of kernels that
$$
\left.\frac{\partial^2P} {\partial
z_b^\beta\partial(z_a^\alpha)^*}\right|_{\mathbf{z}=0}=
\mathrm{const}(q,n)\cdot \left(\frac{1-q^{-2n}}{1-q^{-2}}q^{2(2n-a-\alpha)}
\zeta_a^\alpha(\zeta_b^\beta)^*-\delta_{ab}\delta^{\alpha\beta}\right).
$$
Set here $a=b=c$ to get
$$
\left.\frac{\partial^2P} {\partial
z_c^\beta\partial(z_c^\alpha)^*}\right|_{\mathbf{z}=0}=
\mathrm{const}(q,n)\cdot \left(\frac{1-q^{-2n}}{1-q^{-2}}q^{2(2n-c-\alpha)}\
\zeta_c^\alpha(\zeta_c^\beta)^*-\delta^{\alpha\beta}\right).
$$
On the other hand, the generators of the function algebra on the Shilov
boundary are subject to the relation
$$
\sum\limits_{c=1}^n\zeta_c^\alpha(\zeta^\beta_c)^*=
q^{-2n+\alpha+\beta}\delta^{\alpha\beta},\qquad\alpha,\beta=1,2,\ldots,n,
$$
see \eqref{Shilov_eq}. Hence
\begin{multline*}
\frac1{\mathrm{const}(q,n)}\left.\sum_{c=1}^nq^{2c}\frac{\partial^2P}{\partial
z_c^\beta\partial(z_c^\alpha)^*}\right|_{\mathbf{z}=0}=
\\ =\frac{1-q^{-2n}}{1-q^{-2}}q^{2(2n-\alpha)}
\sum_{c=1}^n\zeta_c^\alpha(\zeta_c^\beta)^*-
q^2\frac{1-q^{2n}}{1-q^2}\delta^{\alpha\beta}=
\\ =\frac{1-q^{-2n}}{1-q^{-2}}q^{2(2n-\alpha)}
q^{-2n+\alpha+\beta}\delta^{\alpha\beta}-q^2
\frac{1-q^{2n}}{1-q^2}\delta^{\alpha\beta}=0,
\end{multline*}
so that the following statement is valid.

\begin{lemma}\label{quantum_Hua_zero}
If $u\in\mathscr{D}(\mathbb{D})_q'$ is a Poisson integral on the quantum
$n\times n$-matrix ball, then
\begin{equation*}
\left.\sum_{c=1}^nq^{2c}\frac{\partial^2u}{\partial
z_c^\beta\partial(z_c^\alpha)^*}\right|_{\mathbf{z}=0}=0
\end{equation*}
for all $\alpha,\beta=1,2,\ldots,n$.
\end{lemma}

Since the subspace of Poisson integrals
$$
u=\int\limits_{S(\mathbb{D})_q}P(\mathbf{z},\zeta)f(\zeta)d\nu(\zeta),\qquad
f\in\mathbb{C}[S(\mathbb{D})]_q,
$$
is a $U_q\mathfrak{sl}_{2n}$-submodule, the above lemma implies Theorem
\ref{quantum_Hua_new}.

It is known \cite{Helg94}, that in the classical case $q=1$ the Poisson kernel $P$ is a
solution of one more equation system
\begin{equation*}
\left.\sum_{\gamma=1}^n\frac{\partial^2u(g\cdot\,\mathbf{z})}{\partial
z_b^\gamma\partial\overline{z}_a^\gamma}\right|_{\mathbf{z}=0}=0,\qquad
g\in SU_{n,n},\quad a,b\in\{1,2,\ldots,n\}.
\end{equation*}
An argument similar to the above allows one to obtain a $q$-analog of this
result.

\begin{proposition}\label{quantum_Hua_new_2}
If $u\in\mathscr{D}(\mathbb{D})_q'$ is a Poisson integral on the quantum
$n\times n$-matrix ball, then
\begin{equation*}
\left.\sum_{\gamma=1}^nq^{2\gamma}\frac{\partial^2(\xi u)}{\partial
z_a^\gamma\partial(z_b^\gamma)^*}\right|_{\mathbf{z}=0}=0
\end{equation*}
for all $\xi\in U_q\mathfrak{sl}_{2n}$, $a,b=1,2,\ldots,n$.
\end{proposition}

\bigskip

\section{Addendum. Hint to a more general case}

Turn from the special case of $n\times n$-matrix ball to a more general
case bounded symmetric domain of tube type. We intend to introduce the Hua
operator, which can be used in order to rewrite the Hua equations in a more
habitual form, see \cite[p. 593]{JohnsonKoranyi}.

Let $\mathfrak{g}$ be a simple complex Lie algebra,
$(a_{ij})_{i,j=1,\ldots,l}$ the associated Cartan matrix. We refer to the
well known (see \cite{Jant}) description of the universal enveloping
algebra $U\mathfrak{g}$ in terms of its generators $e_i$, $f_i$, $h_i$,
$i=1,...,l$, and the standard relations. Consider also the linear span
$\mathfrak{h}$ of the set $\{h_i,i=1,...,l\}$ (a Cartan subalgebra), and
the simple roots $\{\alpha_i\in\mathfrak{h}^*|i=1,...,l\}$ given by
$\alpha_i(h_j)=a_{ji}$. Let $\delta$ be the maximal root,
$\delta=\sum\limits_{i=1}^lc_i\alpha_i$. Assume that it is possible to
choose $l_0\in\{1,...,l\}$ so that $c_{l_0}=1$. Fix an element
$h_{0}\in\mathfrak{h}$ with the following properties:
$$\alpha_i(h_0)=0,\quad i\ne l_0;\qquad\alpha_{l_0}(h_0)=2.$$
In this case the Lie algebra $\mathfrak{g}$ is equipped with the
$\mathbb{Z}$-grading as follows:
\begin{equation}\label{par_type}
\mathfrak{g}=\mathfrak{g}_{-1}\oplus\mathfrak{g}_0\oplus\mathfrak{g}_{+1},
\qquad\mathfrak{g}_j=\{\xi\in\mathfrak{g}|\:[h_0,\xi]=2j\xi\}
\end{equation}
(that is, $\mathfrak{g}_i=\{0\}$ for all $i$ with $|i|>1$).

Denote by $\mathfrak{k}\subset\mathfrak{g}$ the Lie subalgebra, generated
by
$$e_i,f_i,\quad i\ne l_0;\qquad h_i,\quad i=1,...,l.$$
If \eqref{par_type} is true then $\mathfrak{g}_0=\mathfrak{k}$, and the
pair $(\mathfrak{g},\mathfrak{k})$ is called the Hermitian symmetric pair.
In what follows, we obey the conventions of the theory of Hermitian
symmetric spaces, where it is custom to use the notation
$\mathfrak{p}^{\pm}$ instead of $\mathfrak{g}_{\pm 1}$ as in
\eqref{par_type}.

Harish-Chandra introduced a standard realization of an irreducible bounded
symmetric domain $\mathbb{D}$, considered up to biholomorphic isomorphisms,
as a unit ball in the normed space $\mathfrak{p}^-$ \cite{HCh,Wolf_fine}.
Let $G$ be the simply connected complex linear algebraic group with
$\operatorname{Lie}(G)=\mathfrak{g}$, and $K\subset G$ such connected
linear algebraic subgroup that $\operatorname{Lie}(K)=\mathfrak{k}$. In
this context one has the well known Harish-Chandra embedding
$$i:K\backslash G\hookrightarrow\mathfrak{p}^-.$$

Let $W$ be the Weyl group of the root system $R$ of $\mathfrak{g}$, and
$w_0\in W$ the longest element. The irreducible bounded symmetric domain
$\mathbb{D}$ associated to the pair $(\mathfrak{g},\mathfrak{k})$, is a
tube type domain if and only if $\varpi_{l_0}=-w_0\varpi_{l_0}$.

Let $U_q\mathfrak{g}$ be the quantum universal enveloping algebra of
$\mathfrak{g}$. Recall that it is a Hopf algebra and can be described in
terms of its generators $E_i$, $F_i$, $K_i^{\pm 1}$, $i=1,...,l$, and the
standard Drinfeld-Jimbo relations.

Introduce quantum analogs for invariant differential operators to be used
to produce the Hua operator.

Let $V$ be a finite dimensional weight $U_q\mathfrak{g}$-module. In an
obvious way, $V\otimes_{U_q\mathfrak{k}}U_q\mathfrak{g}$ is equipped with a
structure of right $U_q\mathfrak{g}$-module. It is easy to see that in the
category of left $U_q\mathfrak{g}$-modules
\begin{align*}
& \operatorname{Hom}_{U_q\mathfrak{k}}(U_q\mathfrak{g},V)\cong
(V^*\otimes_{U_q\mathfrak{k}}U_q\mathfrak{g})^*, \\ & f \mapsto
\widetilde{f},\qquad\widetilde{f}(l \otimes \xi)=l(f(\xi)),\qquad
l \in V^*,\quad\xi \in U_q \mathfrak{g}.
\end{align*}
The vector space
$\operatorname{Hom}_{U_q\mathfrak{k}}(U_q\mathfrak{g},V)$ is a
quantum analog for the space of sections of a homogeneous vector
bundle on the homogeneous space $K\backslash G$. Suppose we are given
two finite dimensional weight $U_q\mathfrak{g}$-modules $V_1$, $V_2$
and a morphism of right $U_q \mathfrak{g}$-modules
$$
A:V_2^*\otimes_{U_q\mathfrak{k}}U_q\mathfrak{g}\to
V_1^*\otimes_{U_q\mathfrak{k}}U_q\mathfrak{g}.
$$
Associate to the latter morphism the adjoint linear map
$$
A^*:\operatorname{Hom}_{U_q\mathfrak{k}}(U_q\mathfrak{g},V_1)\to
\operatorname{Hom}_{U_q\mathfrak{k}}(U_q\mathfrak{g},V_2),
$$
which is also a morphism of $U_q\mathfrak{g}$-modules. Such dual
operators are treated as quantum analogs of invariant differential
operators.

Thus the invariant differential operators are in one-to-one correspondence
with the elements of the space
\begin{align*}
& \operatorname{Hom}_{U_q\mathfrak{g}}(V_2^*\otimes_{U_q\mathfrak{k}}
U_q\mathfrak{g},V_1^*\otimes_{U_q\mathfrak{k}}U_q\mathfrak{g})\cong
\mathrm{Hom}_{U_q\mathfrak{k}}(V_2^*,V_1^*\otimes_{U_q\mathfrak{k}}
U_q\mathfrak{g}),
\\& f \mapsto \widetilde{f}, \qquad \widetilde{f}(l)=f(l \otimes 1),
\qquad l \in V_2^*.
\end{align*}

Turn to a construction of the Hua operator. Set
$\mathfrak{p}^+=U_q\mathfrak{k}E_{l_0}$,
$\mathfrak{p}^-=U_q\mathfrak{k}(K_{l_0}F_{l_0})$, both are finite
dimensional weight $U_q\mathfrak{k}$-modules \cite{Jak-Hermit}. The
morphisms of right $U_q\mathfrak{k}$-modules
\begin{align*}
\mathfrak{p}^+ &\to\mathbb{C}\otimes_{U_q\mathfrak{k}}U_q\mathfrak{g}, &
E_{l_0} &\mapsto 1\otimes E_{l_0},
\\ \mathfrak{p}^- &\to\mathbb{C}\otimes_{U_q\mathfrak{k}}U_q\mathfrak{g}, &
K_{l_0}F_{l_0} &\mapsto 1\otimes K_{l_0}F_{l_0}
\end{align*}
determine invariant linear differential operators
$$
\operatorname{Hom}_{U_q\mathfrak{k}}(U_q\mathfrak{g},\mathbb{C})\to
\operatorname{Hom}_{U_q\mathfrak{k}}(U_q\mathfrak{g},\mathfrak{p}^\pm).
$$

Recall that $U_q\mathfrak{k}$-modules form a tensor category, and that the
comultiplication $\triangle:U_q\mathfrak{g}\to U_q\mathfrak{g}\otimes
U_q\mathfrak{g}$ is a morphism of this category.

Consider the morphisms of $U_q\mathfrak{g}$-modules
\begin{gather}\label{short_seq_op_1}
\operatorname{Hom}_{U_q\mathfrak{k}}(U_q\mathfrak{g},\mathbb{C})\to
\operatorname{Hom}_{U_q\mathfrak{k}}(U_q\mathfrak{g}\otimes
U_q\mathfrak{g},\mathfrak{p}^+\otimes\mathfrak{p}^-),
\\ \label{short_seq_op_2}
\operatorname{Hom}_{U_q\mathfrak{k}}(U_q\mathfrak{g}\otimes
U_q\mathfrak{g},\mathfrak{p}^+\otimes\mathfrak{p}^-)\to
\operatorname{Hom}_{U_q\mathfrak{k}}(U_q\mathfrak{g},\mathfrak{p}^+\otimes
\mathfrak{p}^-).
\end{gather}
Let $\mathfrak{k}_q$ be the finite dimensional weight $U_q\mathfrak{k}$-module with the
same weights and weight multiplicities as the $U\mathfrak{k}$-module $\mathfrak{k}$.
There exists a unique $U_q\mathfrak{k}$-submodule
$\mathcal{H}_q\subset\mathfrak{p}^+\otimes\mathfrak{p}^-$ such that
$(\mathfrak{p}^+\otimes\mathfrak{p}^-)/\mathcal{H}_q\approx\mathfrak{k}_q$ (because a
similar fact is well known in the classical case $q=1$, see \cite[Proposition
4.2]{BenBurDamRat}). Fix an onto morphism
$\mathfrak{p}^+\otimes\mathfrak{p}^-\to\mathfrak{k}_q$ and consider the associated
invariant `formal' differential operator
\begin{equation}\label{short_seq_op}
\operatorname{Hom}_{U_q\mathfrak{k}}(U_q\mathfrak{g},\mathfrak{p}^+\otimes
\mathfrak{p}^-)\to\operatorname{Hom}_{U_q\mathfrak{k}}
(U_q\mathfrak{g},\mathfrak{k}_q).
\end{equation}
Denote by $\mathcal{D}_q$ the composition of the maps
\eqref{short_seq_op_1}, \eqref{short_seq_op_2}, and
\eqref{short_seq_op}. By definition, $\mathcal{D}_q$ is an invariant
differential operator.

Recall the standard definitions of quantum analogs for the algebras
of regular functions on the group $G$ and on the homogeneous space
$K\backslash G$. Denote by
$\mathbb{C}[G]_q\subset(U_q\mathfrak{g})^*$ the Hopf algebra of all
matrix elements of weight finite dimensional of $U_q\mathfrak{g}$.
$\mathbb{C}[G]_q$ is equipped with a structure of
$U_q^{\mathrm{op}}\mathfrak{g}\otimes U_q\mathfrak{g}$-module algebra
via quantum analogs of the standard right and left regular actions
$(\xi'\otimes\xi'')f=\mathcal{L}_{\mathrm{reg}}(\xi')
\mathcal{R}_{\mathrm{reg}}(\xi'')f$, where
$$
\mathcal{L}_{\mathrm{reg}}(\xi'f)(\eta)=f(\xi'\eta),\qquad
\mathcal{R}_{\mathrm{reg}}(\xi''f)(\eta)=f(\eta\xi''),\qquad\xi',\xi'',\eta
\in U_q\mathfrak{g},\quad f\in\mathbb{C}[G]_q.
$$
($U_q^{\mathrm{op}}\mathfrak{g}$ is the Hopf algebra with the opposite
multiplication). $\mathbb{C}[G]_q$ is called the algebra of regular
functions on the quantum group $G$.

Introduce the notation
$$
\mathbb{C}[K\backslash G]_q=\{\xi\in\mathbb{C}[G]_q|\:
\mathcal{L}_{\mathrm{reg}}(\eta)\xi=0,\quad\eta\in U_q\mathfrak{k}\}.
$$
This Hopf subalgebra is a quantum analog for the algebra of regular
functions on the homogeneous space $K\backslash G$. It is easy to prove
that $\mathbb{C}[K\backslash G]_q\subset
\operatorname{Hom}_{U_q\mathfrak{k}}(U_q\mathfrak{g},\mathbb{C})$, so one
can consider the restriction of $\mathcal{D}_q$ onto
$\mathbb{C}[K\backslash G]_q$.

One can consider a localization $\mathbb{C}[K\backslash G]_{q,x}$ of
the algebra $\mathbb{C}[K\backslash G]_q$ with respect to the Ore set
$x^{\mathbb Z_+}$. It can be proved that the extension of
$\mathcal{D}_q$ up to $\mathbb{C}[K\backslash G]_{q,x}$ is well
defined. Pass from $\mathbb{C}[K\backslash G]_{q,x}$ to
$\operatorname{Pol}(\mathfrak{p}^-)_q$ (via the Harish-Chandra
embedding $\mathcal{I}$, see Section \ref{to-homogeneous} for the
special case) and $\mathscr{D}(\mathbb{D})'_q$ to get a $q$-analog
for the Hua operator.

\end{document}